\documentclass{aic}
\usepackage{pudlatex}
\usepackage{multicol}
\usepackage{amsmath}
\usepackage{hyperref} % nutno texovat 2 x
\newcommand\tw{\mathop{\rm tw}\nolimits}

\usepackage{setspace}

\aicAUTHORdetails{%
  title = {A Lower Bound on the Ramsey Number {\large $R_k(k+1,k+1)$}}, %% please capitalize all significant words
  author = {Pavel Pudl\'ak, Vojt\v{e}ch R\"odl, and William J. Wesley},
  plaintextauthor = {Pavel Pudl\'ak, Vojt\v{e}ch R\"odl, William J. Wesley},
  plaintexttitle = {A Lower Bound on the Ramsey Number R_k(k+1,k+1)}, %%  title without math or LaTeX
}   %%% END \aicAUTHORdetails

\aicEDITORdetails{%
   year={2026},
   number={3},
   received={24 December 2024},   % received date: example: 7 January 2017
   published={24 April 2026},  % published date
   doi={10.19086/aic.2026.3},      % XXX = number of paper, e.g. aic006 for paper#6
}   %%% END \aicEDITORdetails

\begin{document}

\begin{frontmatter}[classification=text]
%% EDITOR: this will force the keywords to appear right after the Abstract.
%%   If the abstract is too long and would force the keywords off the
%%   front page, please comment out % [classification=text] above
%%   This way the keywords will be floated on the bottom of the first page
%%   even though the Abstract spills over to the next page.

%%% AUTHOR: Title goes here.  This line is optional.  You must use it
%%   if title has footnote attached or requires nontrivial typesetting,
%%   e.g., inclusion of linebreaks to force nice layout.
\title{A Lower Bound on the Ramsey Number {\Huge $R_k(k+1,k+1)$}} %% please capitalize all significant words

%%% AUTHOR:
%%% List all authors. If you wish, place grant acknowledgements in \thanks.
%%% In brackets include a short tag for each author.
\author[pudl]{Pavel Pudl\'ak\thanks{Partially supported by GA\v{C}R grant 25-16311S and
    institute's grant RVO~67985840.}}
\author[rodl]{Vojt\v{e}ch R\"odl\thanks{Partially supported by NSF grant DMS 2300347.}}
\author[wesl]{William J. Wesley}

%%% AUTHOR: Abstract goes here
\begin{abstract}
We will prove that $R_k(k+1,k+1)\geq 4 \tw_{\lfloor k/4\rfloor -3}(2)$, where $\tw$ is the tower function defined by ${\tw}_1(x)=x$ and ${\tw}_{i+1}(x)=2^{{\tw}_i(x)}$. We also give proofs of $R_k(k+1,k+2)\geq 4 \tw_{k-7}(2)$, $R_k(k+1,2k+1)\geq 4 \tw_{k-3}(2)$, and $R_k(k+2,k+2)\geq 4 \tw_{k-4}(2)$.
\end{abstract}

\end{frontmatter}

\section{Introduction}
The Ramsey number $R_k(l,m)$ is the minimum $N$ with the property that any colorings of $k$-tuples of $\{1,2\dts N\}$ by red and blue yields either an $l$-set with all its $k$-subsets colored red or an $m$-set with all its $k$-subsets colored blue.

Classic results of Erd\H{o}s and Szekeres \cite{erdos-szekeres} imply that
\[
2^{m/2}\leq R_2(m,m)\leq 2^{2m}.
\]
Several improvements have been made and while the best lower bound in~\cite{spencer} is still of order $2^{(1+o(1))m/2}$, the best upper bound proved recently in~\cite{campos-griffiths-morris-sahasrabudhe} is $(4-\epsilon)^m$ for some $\epsilon>0$. 

For $k\geq 3$, another classic result of Erd\H{o}s and Rado~\cite{erdos-rado} asserts
\bel{classic}
\tw_{k-1}(c_1m^2)\leq R_k(m,m)\leq \tw_k(c_2m),
\ee
where $c_1,c_2$ are positive constants (independent of $m$ and~$k$). While the upper bound holds for every $m\geq k+1$, the lower bound was only proved for $m$ sufficiently large w.r.t.~$k$. 
The exponential gap between lower and upper bounds~(\ref{classic}) is a well-known difficult problem in the area. The lower bound in~(\ref{classic}) is proved by induction. Starting with the probabilistic lower bound on $R_3(m,m)$, the induction argument is based on Erd\H os-Hajnal \emph{stepping-up lemma} that allows one to bound $R_{k+1}(m,m)$ in terms of $R_k(m',m')$ for some $m'$ depending on $m$~\cite{erdos-hajnal-rado,GRS}. The stepping-up lemma has been improved by Conlon, Fox, and Sudakov in \cite{CFS}. While the original version required the restriction $m\geq m_0(k)$ with $m_0(k)$ exponential in $k$, the improved version from~\cite{CFS} can be applied as long as $m\geq \frac 52k+4$.\footnote{This bound is not explicitly stated in~\cite{CFS}, but is a direct consequence of their Theorem~4 combined with the base case proved in~\cite{mckay}.} 
For $m$ closer to~$k+1$, the stepping-up lemma so far does not work. Some improvements are also known for off-diagonal Ramsey numbers.
A recent result in this direction is the bound
\[
R_k(k+1,m)\geq \tw_{k-2}(m^{c\log m}),
\]
where $c$ is a constant and $m\geq k!2^k$, proved by Mubayi and Suk~\cite{mubayi-suk}.\footnote{In that paper the bound is stated as $m>k$, but this is a typo because their proof requires an exponential lower bound as we learned from the authors.}

Lower bounds on $k$-hypergraph Ramsey numbers can also be proved using lower bounds on Ramsey numbers for ordered paths in the shift graph on $k$-sets. Let's denote by $P_k(m,n)$ the minimum number $N$ such that any two-coloring of $k$-sets yields an ordered path of length~$m$ in the first color, or an ordered path of length $n$ in the second color in the shift graph of $k$-sets. Clearly,
\bel{P<R}
P_k(m,n)\leq R_k(k+m-1,k+n-1).
\ee
Moshkovitz and Shapira developed a method by which one can compute the Ramsey numbers of ordered paths $P_k(m,n)$. In~\cite{moshkovitz} they considered the diagonal case, which was subsequently generalized to the of-diagonal path Ramsey numbers by Milans, Stolee, and West~\cite{milans}. Their method together with inequality~(\ref{P<R}) enables one to prove lower bounds on the Ramsey numbers. The bounds obtained in this way are smaller than those obtained by the stepping-up lemma, but the method works in the range where the stepping-up lemma fails. However, for $R_k(k+1,k+1)$ this method fails because $P_k(2,2)\leq 2k+1$.

In this paper we will prove that $R_k(k+1,k+1)$ is lower-bounded by a tower of height approximately $k/4$. To this end we will use colorings of shift graphs by three colors. The proof is based on an idea that we used in~\cite{colorings} and a lemma about a certain hypergraph.
%We use a coloring of this hypergraph by two colors. 

Our result is proved in the main part. In Appendix, we show how to prove tower-type lower bounds on the Ramsey numbers $R_k(k+1,k+2), R_k(k+2,k+2)$ and $R_k(k+1,2k+1)$ using a 3-coloring of the shift graph.

%%%%%%%%%%%%%%%%%%%%%%%%%%%%%%%%%%%%%%%%%%%%%%%%%%%%%%%%%%%%%%%%%%%%%%%%%%%%%

\section{Preliminaries}

A shift graph $Sh(N,k)$ is a graph on ${N\choose k}$ vertices corresponding to $k$-element subsets of $[N]$, where, as usual, $[N]=\{1,2\dts N\}$. Two vertices in the shift graph $\{x_1,x_2\dts x_k\}$, $\{y_1,y_2\dts y_k\}$, $x_1<x_2\dots x_k$, $y_1<y_2<\dots y_k$, are joined by an oriented edge if $x_{i+1}=y_i$ for all $i=1,2\dts k-1$. Orientation is important for the Ramsey numbers of paths, but when we talk about the chromatic number, we can forget orientation. 
For a graph $G$, let $\chi(G)$ be its chromatic number.

It follows from the Ramsey theorem that for every $k$, $\chi(Sh(k,N))$ tends to infinity as $N$ goes to infinity. If $N\geq 2k+1$, then $\chi(Sh(k,N))\geq 3$, however, the shift graph remains $3$-colorable for large $N$. The colorings of large shift graphs with a small number of colors has been used for estimating Ramsey numbers and we will use them also here. To this end we introduce the function $s(k)$ by
\[
s(k):=\max\{N;\ \chi(Sh(N,k))\leq 3\}.
\]
It is well-known that $s(k)$ is growing with $k$ roughly as the tower function. In~\cite{colorings}, we proved
\bel{e-bound}
s(k)\geq 4\ \tw_{k-3}(2),
\ee
but our aim in that paper was only to give a simple proof of a tower-function lower bound. Therefore this bound is certainly not optimal that one can prove using known methods. 
Since our bound on Ramsey numbers is based on this function and we do not know the precise value of $s(k)$, we prefer to state our result in terms of this function. 

\bt
For every $k\geq 4$,
$R_k(k+1,k+1)\geq s(\lfloor k/4\rfloor)$.
\et

%%%%%%%%%%%%%%%%%%%%%%%%%%%%%%%%%%%%%%%%%%%%%%%%%%%%%%%%%%%%%%%%%%%%%%%%

\section{The lower bound $R_k(k+1,k+1)\geq s(\lfloor k/4\rfloor)$}

The proof of this inequality uses an idea from our paper \cite{colorings} where we constructed colorings with low discrepancy. Our aim here is different; we only want to prove a lower bound.

Let $k\geq 4$. Because of the monotonicity of $R_k(k+1,k+1)$, we can assume that $k$ is divisible by~4. Let $N:=s(k/4)$. By the definition of~$s$, this means that the chromatic number of the shift graph $Sh(N,k/4)\leq 3$ (for $k$ sufficiently large equality holds true). Our coloring of ${[N]\choose k}$ will be a composition of two mappings.

The first mapping is defined as follows. Let $\phi$ be a 3-coloring of $Sh(N,k/4)$. Divide each $X\in{[N]\choose k}$ into consecutive segments of the same size $l:=k/4$, which we denote by $X^1,X^2,X^3,X^4$. For each  $X\in{[N]\choose k}$, we define
\[
\lambda(X):=(\phi(X^1),\phi(X^2),\phi(X^3),\phi(X^4)).
\]
In this way, we have defined a mapping $\lambda:{[N]\choose k}\to [3]^4$ assigning a vector of length~4 to every~$X$.

Our next goal is to find a 2-coloring of $[3]^4$, i.e., $\psi:[3]^4\to[2]$, such that $\psi\circ\lambda:{[N]\choose k}\to[2]$ has no monochromatic copy of $K^{(k)}_{k+1}$. We will prove a stronger statement:
\begin{quote}
\emph{For every $k+1$ element subset $Y$ of $[N]$, there are five special $k$-subsets of $Y$ on which $\psi\circ\lambda$ is not monochromatic.}
\end{quote}
Let $Y:=\{y_1<y_2<\dots <y_{k+1}\}$ and $k=4l+1$; the five sets are
\bel{e-yz}
Z_0=Y\setminus\{y_1\},\ Z_1=Y\setminus\{y_{l+1}\},\ Z_2=Y\setminus\{y_{2l+1}\},\ Z_3=Y\setminus\{y_{3l+1}\},\ Z_4=Y\setminus\{y_{4l+1}\}.
\ee
On such 5-tuples, the vectors $(\lambda(Z_0),\lambda(Z_1),\lambda(Z_2),\lambda(Z_3),\lambda(Z_4))$ have a special structure, which we call a \emph{bridge}.

\bdf[bridge]
For any $\bar a=(a_1,a_2,a_3,a_4)\in[3]^4$ and $b=(b_1,b_2,b_3,b_4)\in[3]^4$ with $a_i\neq b_i$ for $i=1,2,3,4$, the set of vectors %containing $\bar a$, $\bar b$ and the vectors
\[\begin{array}{l}
(a_1,a_2,a_3,a_4),\\(a_1,a_2,a_3,b_4),\\(a_1,a_2,b_3,b_4),\\(a_1,b_2,b_3,b_4),\\(b_1,b_2,b_3,b_4)
\end{array}\]
is called a \emph{bridge}.
\edf

\bll{l-bridge}
For every $Y\in{[n]\choose k+1}$ and $Z_0\dts Z_4$ as defined in~(\ref{e-yz}), the set of vectors $\lambda(Z_0)\dts\lambda(Z_4)$ forms a bridge.
\el
\bprf
Recalling that $Y=\{y_1<y_2\dots <y_{4l+1}\}$, for $j=0,1,2,3$, set 
\[
A_j:=\{y_{jl+2},\dots,y_{(j+1)l+1}\},\quad B_j:=\{y_{jl+1},\dots,y_{(j+1)l}\}
\]
Observe that, for $i=0,1,2,3,4$, $Z_i$ is the concatenation of the first $i$ blocks $B_j$ followed by the last $4-i$ blocks $A_{j'}$. In particular $Z_0=A_0A_1A_2A_3$ and $Z_4=B_0B_1B_2B_3$. Also, for each~$j$, $0\leq j\leq 3$, the blocks $A_j$ and $B_j$ are in a shift position and hence $\phi(A_j)\neq\phi(B_j)$. Consequently, $\lambda(Z_0),\lambda(Z_1),\dots,\lambda(Z_4)$ form a bridge.
\eprf
The following is a key lemma.
\bll{l-key}
There exists a 2-coloring of $[3]^4$ such that no bridge is monochromatic.
\el
%% We found the coloring by computer search. At present, we do not have a concise definition of such a coloring nor an explanation of why it exists. The coloring is in the Appendix 2, and we will explain how we found it in Section~\ref{how}.
%%%%%%%%%%%%%%%%%%%%%%%%%%%%%%%%%%%%%%%%%%%%%%%%%%%%%%%%%%%%%%%%%%%%%%%%%%%%%%
\bprf
We will show that the coloring $\chi : \Z_3^4 \to \{1,2\}$ given by 
\[
    \chi(v) = \chi((v_1,v_2,v_3,v_4))=  \begin{cases}
          1 \mbox{ if } v_1+v_2+v_3+v_4 = 0 \mbox{ or } v_1+v_3 = 0, \\
          2 \mbox{ otherwise,} 
    \end{cases}
\]
contains no monochromatic bridges.

     Let $\{(a_1,a_2,a_3,a_4),(a_1,a_2,a_3,b_4),(a_1,a_2,b_3,b_4),(a_1,b_2,b_3,b_4),(b_1,b_2,b_3,b_4) \}$ be a bridge. Let $a = (a_1,a_2,a_3,a_4), b = (b_1,b_2,b_3,b_4)$. We claim this bridge is not monochromatic. 

    Case 1: Suppose $\chi(a) = 1$ and $a_1 + a_2 +a_3 + a_4 = 0$. Then since $a_4 \neq b_4$, we have $a_1+a_2+a_3+b_4 \neq 0$. So either $\chi((a_1,a_2,a_3,b_4)) = 2$ and the bridge is not monochromatic, or $a_1 + a_3 = 0$. In the latter case, since $a_3 \neq b_3$, we have $a_1 + b_3 \neq 0$. Again, either the bridge is not monochromatic, or $a_1+a_2 + b_3 + b_4 = 0$. But since $a_2 \neq b_2$, we have that $a_1+b_2+b_3+b_4 \neq 0$, so $\chi((a_1,b_2,b_3,b_4)) = 2$ and the bridge is not monochromatic. 

    Case 2: Suppose $\chi(a) = 1$ and $a_1 +a_3 = 0$. Then $a_1 + b_3 \neq 0$, so if the bridge is monochromatic, we must have $a_1 + a_2 + b_3 + b_4 = a_1+b_2+b_3+b_4  = 0$, but this is impossible since $a_2 \neq b_2$, so by contradiction, the bridge is not monochromatic. 

    Case 3: Suppose $\chi(a) = 2$ and the bridge is monochromatic. Then for each element in the bridge, the sum of the coordinates cannot be zero, and since $a_i \neq b_i$ for all $i$, the sums of the coordinates of ``consecutive" bridge elements alternate values. That is, for some $c \in \{-1,1\}$, we have
 \begin{align}
        a_1+a_2+a_3+a_4 &= c\\ 
        a_1+a_2+a_3+b_4 &= -c \\
        a_1+a_2+b_3+b_4 &= c\\ 
        a_1+b_2+b_3+b_4 &= -c \\
        b_1+b_2+b_3+b_4 &= c.
    \end{align}
    We have also that $a_1 + a_3 \neq 0, a_1+b_3 \neq0$, and $b_1+b_3 \neq 0$. It must be that for some $d \in \{-1,1\}$ we have
    \begin{align}
        a_1+a_3 &= d \\
        a_1+b_3 &= -d \\
        b_1+b_3 &= d,
    \end{align} because otherwise we would have $a_3 = b_3$ or $a_1 = b_1$. 

    Equations $(6)$ and $(7)$ imply $b_3 - a_3 = 2c $, and equations $(10)$ and $(11)$ imply $b_3-a_3 = -2d$. Combining these gives $c = -d$. Now subtracting equation $(10)$ from equation $(6)$ gives $a_2 + a_4 = -c - d =0$, and subtracting equation $(11)$ from equation $(7)$ gives $a_2+b_4 = c + d =0$. Thus $a_4 = b_4$, which contradicts that $a_4$ and $b_4$ are distinct.
\end{proof}

The bound on $R_k(k+1,k+1)$ now follows from Lemmas~\ref{l-bridge} and~\ref{l-key}.

%%%%%%%%%%%%%%%%%%%%%%%%%%%%%%%%%%%%%%%%%%%%%%%%%%%%%%%%%%%%%%%%%%%%%

\section{The way to the lower bound on {\Large $R_k(k+1,k+1)$.}}\label{how}

Consider the following generalization of the concept of a bridge.
\bdf[$n$-bridge]
Let $n\geq 2$. 
For any $\bar a=(a_1\dts a_n)\in[3]^n$ and $\bar b=(b_1\dts b_n)\in[3]^n$ with $a_i\neq b_i$ for $i=1\dts n$, the set of vectors 
\[\begin{array}{l}
(a_1,a_2\dts a_{n-1},a_n),\\ (a_1,a_2\dts a_{n-1},b_n),\\ (a_1,a_2\dts b_{n-1},b_n),\\ \dots\\ (a_1,b_2\dts b_{n-1},b_n),\\ (b_1,b_2\dts b_{n-1},b_n) 
\end{array}\]
is called an $n$-\emph{bridge}.
\edf
Denote by ${\cal B}_n$ the hypergraph on $[3]^n$ whose edges are $n$-bridges. When trying to prove a lower bound on $R_k(k+1,k+1)$, Pudl\'ak and R\"odl first realized the following fact:
\begin{claim}
If ${\cal B}_n$ is 2-colorable, then $R_k(k+1,k+1)\geq s(\lfloor k/n\rfloor)$.
\end{claim}
It turned out that $n=4$ is the least $n\geq 2$ such that ${\cal B}_n$ is $2$-colorable. It is easy to prove that ${\cal B}_2$ is not $2$-colorable; for proving that ${\cal B}_3$ is not 2-colorable and ${\cal B}_4$ is colorable, they used a computer. Specifically, they transformed the problems into problems about satisfiability of CNF formulas and used the SAT-solver \emph{MiniSat,}
\cite{minisat}. About a year later, Wesley found the 2-coloring of ${\cal B}_4$ presented above and a proof that it is correct. Thus he eliminated the need of using a computer to verify the result. Interestingly, he also used a SAT-solver to find the coloring.

%% To represent the problem of coloring a hypergraph $H=(V,E)$ by two colors as a problem about satisfiability of a CNF formula is easy: For every $v\in V$ take a variable $x_v$ and for every hyperedge $\{v_1\dts v_r\}$ take the two clauses $\{x_{v_1}\dts x_{v_r}\},\{\bar x_{v_1}\dts \bar x_{v_r}\}$.

One could try to get a better bound by using shift graphs that are colored by more than 3~colors. For instance, let $s_4(k)$ be the maximal $N$ such that $\chi(Sh(N,k))\leq 4$. One can show a larger bound on $s_4$ than we have for~$s$; namely,
$s_4(k)\geq 4\tw_{k}(2).$
Let ${\cal B}^4_4$ be the natural extension of ${\cal B}_4$ to a hypergraph on $[4]^4$. If ${\cal B}^4_4$ were 2-colorable, we would get a lower bound $\geq s_4(\lfloor k/4\rfloor)$. But this hypergraph is not 2-colorable, as we found out using the SAT-solver. So our bound is the best possible that one can get using shift graphs and bridges.

\section{Conclusions}

We have made progress in computing the value of $R_k(k+1,k+1)$, but our estimate is certainly far from the true value. We believe that $R_k(k+1,k+1)$ is at least the tower of $k-c$ twos for some constant $c$. A possible reason why we cannot get such a bound using our method is that we only consider a constant number of $k$-sets in a potentially homogeneous $(k+1)$-set~$X$. Specifically our coloring is such that there is a configuration of five particular $k$-subsets of $X$ which is not monochromatic for any $X$, $|X|=k+1$.\footnote{But note that for the lower bound on $R_k(k+2,k+2)$ we only needed \emph{three} particular sets to get almost full tower of twos.}

%%% AUTHOR: optional appendix here
\appendix %% you may comment this out if no Appendix
\section*{Appendix}
\setcounter{section}{1}
%\section{Lower bounds on paths and Ramsey numbers}

In this appendix, we will prove lower bounds on $R_k(k+1,k+2)$, $R_k(k+2,k+2)$ and $R_k(k+1,2k+1)$. Our lower bounds are based on colorings of shift graphs.
We prove the lower bounds on $R_k(k+1,k+2)$, $R_k(k+2,k+2)$ by proving lower bounds on the Ramsey number of paths $P_k(2,3)$ and $P_k(3,3)$. As we learned from Dhruv Mubayi, similar bounds can be computed using the method of Moshkovitz and Shapira. We include these bounds because our proofs based on 3-coloring of shift graphs are very simple. These bounds can also be viewed as showing relations between the numbers  $P_k(2,3)$ and $P_k(3,3)$ on the one hand, and $P_l(2,2,2)$ on the other, since $P_l(2,2,2)=s(l)+1$. The bound on $R_k(k+2,k+2)$ appeared in~\cite{duffus-lefmann-rodl}; here we give a simplified proof.

\bpr
For every $k\geq 4$,
\ben
\item $R_k(k+1,k+2)\geq P_k(2,3)\geq s(k-4)$,
\item $R_k(k+2,k+2)\geq P_k(3,3)\geq s(k-1)$,
\item $R_k(k+1,2k+1)\geq s(k)$.
\een
\epr

\paragraph{Proof of $P_k(2,3)\geq s(k-4)$.}
Let $N:=s(k-4)$ and let $\phi:{[N]\choose k-4}\to\{1,2,3\}$ be a 3-coloring of the shift graph $Sh(N,k-4)$. We color $k$-sets with two colors as follows.
Let $X=\{x_1,x_2\dts x_k\}$, and consider the colors of the five $(k-4)$-sets
\begin{align*}
& c_1=\phi(\{x_1,x_2\dts x_{k-4}\}),\\
& c_2=\phi(\{x_2,x_3\dts x_{k-3}\}),\\
& \dots\qquad\\
& c_5=\phi(\{x_5,x_6\dts x_{k}\}).
\end{align*}
We color $X$ red, if
\ben
\item either $c_2<c_3>c_4$, (type I.),
\item or $c_1>c_2>c_3<c_4<c_5$, (type II.).
\een
Otherwise $X$ is blue.

First, we show that there is no red path of length~2. It is clear that such a path cannot have the two $k$-sets of the same type. Suppose that $X$ and $Y$ are red, they are in a shift position, $X$ is before $Y$, and $X$ is type I and $Y$ is type~II. Then the middle three $(k-4)$-sets in $X$ become the initial three $(k-4)$-sets in~$Y$, which is impossible, because the colors of the middle three sets of $X$ are not monotonic, while in $Y$ they monotonically decrease. The case of $X$ being type II and $Y$ type~I is essentially the same.

Now suppose there is a blue path $X_1,X_2,X_3$. Consider the path of length seven of $(k-4)$-subsets in $X_1\cup X_2\cup X_3$. Let $c_1,c_2,\dots,c_7$ be the colors of the consecutive sets in this path. Recall that $c_i\neq c_{i+1}$ for $i=1,2\dts 6$ because consecutive $(k-4)$ sets are in shift positions. First, note that none of the colors $c_3,c_4,c_5$ can be equal to~3 because otherwise $X_1$, or $X_2$, or $X_3$ would be red, type~I. Further, $c_4$ cannot equal to~2; this is because, since neither $c_3$ nor $c_5$ equals to~$3$, $c_4=2$ would imply that $c_3=c_5=1$ contradicting the fact that $X_2$ is not red, type~I. Having established $c_4\not\in\{2,3\}$ and $c_3\neq 3$ and $c_5\neq 3$, we infer that  $c_4=1$ and hence $c_3=c_5=2$.

Next observe that $c_2\neq 1$, for otherwise $X_1$ would be red of type I (with $c_2=1,c_3=2,c_4=1$). Consequently, since $c_2\neq c_3=2$, $c_2$ must be equal to~3. For the same reason, we have also $c_6=3$. But then $c_2>c_3>c_4<c_5<c_6$ contradicting the assumption that $X_2$ is not red, type~II.
\qed

\paragraph{Proof of $P_k(3,3)\geq s(k-1)$.}

Let $N:=s(k-1)$ and let $\phi:{[N]\choose k-1}\to\{1,2,3\}$ be a 3-coloring of the shift graph $Sh(N,k-1)$. We will define a coloring of $k$-sets by two colors, red and blue, as follows. For a $k$-set $X=\{x_1<\dots <x_k\}$, we color $X$ red if $\phi\{x_1\dts x_{k-1}\}<\phi\{x_2\dts x_{k}\}$, and blue if $\phi\{x_1\dts x_{k-1}\}>\phi\{x_2\dts x_{k}\}$. Equality is not possible because $\phi$ is a proper coloring of the shift graph. Let $Y:=\{y_1<\dots<y_{k+2}\}$ be a ($k+2$)-subset of~$[N]$ and assume that 
 $\{y_1\dts y_k\},\{y_2\dts y_{k+1}\},\{y_3\dts y_{k+2}\}$ are colored red. Then we have
\bel{e-increasing}
\phi\{y_1\dts y_{k-1}\}<\phi\{y_2\dts y_{k}\}<\phi\{y_3\dts y_{k+1}\}<\phi\{y_4\dts y_{k+2}\}.
\ee
% because the sets $\{y_1\dts y_k\},\{y_2\dts y_{k+1}\},\{y_3\dts y_{k+2}\}$ are colored red.
But this is impossible, because $\phi$ has only 3~colors. The same argument shows that  $$\{y_1\dts y_k\},\{y_2\dts y_{k+1}\},\{y_3\dts y_{k+2}\}$$ cannot be all colored blue.
\qed

\paragraph{Proof of  $R_k(k+1,2k+1)\geq s(k)$.}

Let $N:=s(k)$ and let $\phi:{[N]\choose k}\to\{1,2,3\}$ be a 3-coloring of the shift graph $Sh(N,k)$. Our coloring of $k$-sets by two colors, red and blue, is defined from $\phi$ in a very simple way: we color a $k$-set $X$ red if $\phi(X)=1$, otherwise it is blue. Since no pair of $k$-sets in a shift position has the same color by $\phi$, no $K^{(k)}_{k+1}$ is colored red. On the other hand, the subgraph of $Sh(N,k)$ induced by blue vertices is bipartite (such a 2-coloring is~$\phi$). Since every $(2k+1)$-subset $Y$ of $[N]$ induces an odd cycle in the shift graph, ${Y\choose k}$ cannot be all colored blue.
\qed

%%% AUTHOR: optional acknowledgments here
\section*{Acknowledgments} %%  you may comment this out if no Ackno
We thank Sam Buss and Neil Thapen for helping us with SAT solving, and Petr Pudl\'ak for writing the code for generating clauses. We also thank Dhruv Mubayi for informing us about the results on Ramsey numbers of paths.

%%%%%%%%%%%%%%%%%%%%%%%%%%%%%%%%%%%%%%%%%%%%%%%%%%%%%%%%%%%%%%%%%%%%%%%%%%%%

%\newpage
\begin{aicauthors}
\begin{authorinfo}[pudl]
  Pavel Pudl\'ak\\
  Institute of Mathematics, CAS\\
  Prague, Czech Republic\\
  pudlak\imageat{}math\imagedot{}cas\imagedot{}cz \\
  \url{https://users.math.cas.cz/~pudlak}
\end{authorinfo}
\begin{authorinfo}[rodl]
  Vojt\v{e}ch R\"odl\\
  Emory University\\
  Atlanta, USA\\
  vrodl\imageat{}emory\imagedot{}edu \\
  \url{https://www.math.emory.edu/~rodl}
\end{authorinfo}
\begin{authorinfo}[wesl]
  William J. Wesley\\
  University of California, San Diego\\
  La Jolla, USA\\
  wjwesley\imageat{}ucsd\imagedot{}edu\\
  \url{https://mathweb.ucsd.edu/~wjwesley}
\end{authorinfo}
\end{aicauthors}

\end{document}